\documentclass[12pt]{article}
\usepackage{a4}
\usepackage{amsthm}
\usepackage{amsfonts}
\usepackage{amssymb}
\usepackage{amsmath}
\usepackage{cite}
\usepackage{epsfig}
\usepackage{tikz}
\usepackage{filecontents}
\usepackage{hyperref}

\newtheorem{theorem}{Theorem}
\newtheorem{corollary}[theorem]{Corollary}
\newtheorem{lemma}[theorem]{Lemma}

\newtheorem{conjecture}[theorem]{Conjecture}
\theoremstyle{definition}

\def\EE{{\mathbb E}\;}

\def\EEW{{\mathbb E}_W}
\def\NN{{\mathbb N}}

\def\FF{{\mathcal F}}

\begin{document}
\title{Decomposing graphs into edges and triangles\thanks{This work has received funding from the European Research Council (ERC) under the European Union’s Horizon 2020 research and innovation programme (grant agreement No 648509). This publication reflects only its authors' view; the European Research Council Executive Agency is not responsible for any use that may be made of the information it contains.}}

\author{Daniel Kr\'al'\thanks{Mathematics Institute, DIMAP and Department of Computer Science, University of Warwick, Coventry CV4 7AL, UK. E-mail: {\tt d.kral@warwick.ac.uk}. The first author was also supported by the Engineering and Physical Sciences Research Council Standard Grant number EP/M025365/1.}\and
        Bernard Lidick\'y\thanks{Department of Mathematics, Iowa State University. Ames, IA, USA. E-mail: {\tt lidicky@iastate.edu}. This author was supported in part by NSF grant DMS-1600390.}\and
        Ta\'isa L.~Martins\thanks{Mathematics Institute, University of Warwick, Coventry CV4 7AL, UK. E-mail: {\tt t.lopes-martins@warwick.ac.uk}. This author was also supported by the CNPq Science Without Borders grant number 200932/2014-4.}\and
 Yanitsa Pehova\thanks{Mathematics Institute, University of Warwick, Coventry CV4 7AL, UK. E-mail: {\tt y.pehova@warwick.ac.uk}.}}
\date{} 
\maketitle

\begin{abstract}
We prove the following 30-year old conjecture of Gy\H ori and Tuza:
the edges of every $n$-vertex graph $G$ can be decomposed into complete graphs $C_1,\ldots,C_\ell$ of orders two and three
such that $|C_1|+\cdots+|C_\ell|\le (1/2+o(1))n^2$.
This result implies the asymptotic version of the old result of Erd\H os, Goodman and P\'osa that
asserts the existence of such a decomposition with $\ell\le n^2/4$.
\end{abstract} 

\begin{flushleft}
\end{flushleft}

\section{Introduction}

Results on the existence of edge-disjoint copies of specific subgraphs in graphs
is a classical theme in extremal graph theory.
Motivated by the following result of Erd\H os, Goodman and P\'osa~\cite{bib-erdos66+},
we study the problem of covering edges of a given graph by edge-disjoint complete graphs.
\begin{theorem}[Erd\H os, Goodman and P\'osa~\cite{bib-erdos66+}]
\label{thm-erdos+}
The edges of every $n$-vertex graph can be decomposed into at most $\lfloor n^2/4\rfloor$ complete graphs.
\end{theorem}
In fact, they proved the following stronger statement.
\begin{theorem}[Erd\H os, Goodman and P\'osa~\cite{bib-erdos66+}]
\label{thm-erdos++}
The edges of every $n$-vertex graph can be decomposed into at most $\lfloor n^2/4\rfloor$ copies of $K_2$ and $K_3$.
\end{theorem}
The bounds given in Theorems~\ref{thm-erdos+} and~\ref{thm-erdos++} are best possible
as witnessed by complete bipartite graphs with parts of equal sizes.

Theorem~\ref{thm-erdos+} actually holds in a stronger form that we now present.
Chung~\cite{bib-chung81}, Gy\H ori and Kostochka~\cite{bib-gyori79+}, and Kahn~\cite{bib-kahn81}, independently,
proved a conjecture of Katona and Tarj\'an asserting that the edges of every $n$-vertex graph
can be covered with complete graphs $C_1,\ldots,C_\ell$ such that the sum of their orders is at most $n^2/2$.
In fact, the first two proofs yield a stronger statement,
which implies Theorem~\ref{thm-erdos+} and which we next state as a separate theorem.
To state the theorem, 
we define $\pi_k(G)$ for a graph $G$ to be the minimum integer $m$ such that the edges of $G$
can be decomposed into complete graphs $C_1,\ldots,C_\ell$ of order at most $k$
with $|C_1|+\cdots+|C_\ell|=m$, and we let $\pi(G)=\min_{k\in\NN}\pi_k(G)$.
\begin{theorem}[Chung~\cite{bib-chung81}; Gy\H ori and Kostochka~\cite{bib-gyori79+}]
\label{thm-all}
Every $n$-vertex graph $G$ satisfies $\pi(G)\le n^2/2$.
\end{theorem}
Observe that Theorem~\ref{thm-all} indeed implies the existence of a decomposition into at most $\lfloor n^2/4\rfloor$ complete graphs.
McGuinnes~\cite{bib-mcguinness94a,bib-mcguinness94b} extended these results
by showing that decompositions from Theorems~\ref{thm-erdos+} and~\ref{thm-all} can be constructed in the greedy way,
which confirmed a conjecture of Winkler of this being the case in the setting of Theorem~\ref{thm-erdos+}.

In view of Theorem~\ref{thm-erdos++},
it is natural to ask whether Theorem~\ref{thm-all} holds under the additional assumption that 
all complete graphs in the decomposition are copies of $K_2$ and $K_3$,
i.e., whether $\pi_3(G)\le n^2/2$.
Gy\H ori and Tuza~\cite{bib-gyori87+} provided a partial answer by proving that $\pi_3(G)\le 9n^2/16$ and
conjectured the following.
\begin{conjecture}[Gy\H ori and Tuza~{\cite[Problem 40]{bib-tuza91}}] %
\label{conj}
Every $n$-vertex graph $G$ satisfies $\pi_3(G)\le (1/2+o(1))n^2$.
\end{conjecture}
We prove this conjecture.
Our result also solves~\cite[Problem 41]{bib-tuza91}, which we state as Corollary~\ref{cor-main}.
We remark that we stated the conjecture in the version given by Gy\H ori in several of his talks and
by Tuza in \cite[Problem 40]{bib-tuza91};
the paper~\cite{bib-gyori87+} contains a version with a different lower order term.

We would also like to mention a closely related variant of the problem suggested by Erd\H os,
where the cliques in the decomposition have weights one less than their orders.
Formally, 
define $\pi^-(G)$ for a graph to be the minimum $m$ such that the edges of a graph $G$
can be decomposed into complete graphs $C_1,\ldots,C_\ell$ with $(|C_1|-1)+\cdots+(|C_\ell|-1)=m$.
Erd\H os asked, see~\cite[Problem 43]{bib-tuza91}, whether $\pi^-(G)\le n^2/4$ for every $n$-vertex graph $G$.
This problem remains open and
was proved for $K_4$-free graphs only recently by Gy\H ori and Keszegh~\cite{bib-gyori15+,bib-gyori17+}.
Namely, they proved that
every $K_4$-free graph with $n$ vertices and $\lfloor n^2/4\rfloor+k$ edges contains $k$ edge-disjoint triangles.

\section{Preliminaries}

We follow the standard graph theory terminology;
we review here some less standard notation and
briefly introduce the flag algebra method.
If $G$ is a graph, then $|G|$ denotes the order of $G$, i.e., the number of vertices of $G$.
Further, if $W$ is a subset of vertices of $G$, then $G[W]$ is the subgraph of $G$ induced by $W$.

In our arguments, we also consider fractional decompositions.
A fractional $k$-decomposition of a graph $G$ is an assignment of non-negative real weights
to complete subgraphs of order at most $k$ such that
the sum of the weights of the complete subgraphs containing any edge $e$ is equal to one.
The weight of a fractional $k$-decomposition is the sum of the weights of the complete subgraphs multiplied by their orders, and
the minimum weight of a fractional $k$-decomposition of a graph $G$ is denoted by $\pi_{k,f}(G)$.
Observe that $\pi_{k,f}(G)\le\pi_k(G)$ for every graph $G$.

\subsection{Flag algebra method}

The flag algebra method introduced by Razborov~\cite{bib-razborov07} has changed
the landscape of extremal combinatorics.
It has been applied to many long-standing open problems,
e.g.\cite{bib-baber11,bib-baber14,bib-balogh17,bib-balogh16,bib-balogh14,bib-balogh15,bib-coregliano17,bib-das13,bib-evan15,bib-falgas15,bib-gethner17,bib-glebov16,bib-goaoc15,bib-hatami12,bib-hladky17,bib-kral12,bib-kral13,bib-kral14,bib-lidicky17+,bib-cumm13}.
The method is designed to analyze asymptotic behavior of substructure densities and
we now briefly describe it.

We start by introducing some necessary notation.
The family of all finite graphs is denoted by $\FF$ and
the family of graphs with $\ell$ vertices by  $\FF_\ell$.
If $F$ and $G$ are two graphs, then $p(F,G)$ is the probability that
$|F|$ distinct vertices chosen uniformly at random among the vertices of $G$ induce a graph isomorphic to $F$;
if $|F|>|G|$, we set $p(F,G) = 0$.
A type is a graph with its vertices labeled with $1,\ldots,|\sigma|$ and
a $\sigma$-flag is a graph with $|\sigma|$ vertices labeled by $1,\ldots,|\sigma|$ such that
the labeled vertices induce a copy of $\sigma$ preserving the vertex labels.
In the analogy with the notation for ordinary graphs,
the set of all $\sigma$-flags is denoted by $\FF^\sigma$ and
the set of all $\sigma$-flags with exactly $\ell$ vertices by $\FF^\sigma_\ell$.

We next extend the definition of $p(F,G)$ to $\sigma$-flags and generalize it to pairs of graphs.
If $F$ and $G$ are two $\sigma$-flags,
then $p(F,G)$ is the probability that 
$|F|-|\sigma|$ distinct vertices chosen uniformly at random among the unlabeled vertices of $G$
induce a copy of the $\sigma$-flag $F$;
if $|F|>|G|$, we again set $p(F,G) = 0$.
Let $F$ and $F'$ be two $\sigma$-flags and $G$ a $\sigma$-flag with at least $|F|+|F'|-|\sigma|$ vertices.
The quantity $p(F,F'; G)$ is the probability that
two disjoint $|F|-|\sigma|$ and $|F'|-|\sigma|$ subsets of unlabeled vertices of $G$
induce together with the labeled vertices of $G$ the $\sigma$-flags $F$ and $F'$, respectively.
It holds~\cite[Lemma 2.3]{bib-razborov07} that
\begin{equation}
p(F,F'; G)  = p(F,G) \cdot p(F',G) + o(1) \label{eq:prodX}
\end{equation}
where $o(1)$ tends to zero with $|G|$ tending to infinity.

Let $\vec F = [F_1,\ldots,F_t]$ be a vector of $\sigma$-flags, i.e., $F_i \in \FF^\sigma$. 
If $M$ is a $t\times t$ positive semidefinite matrix, it follows from \eqref{eq:prodX}, see~\cite{bib-razborov07}, that
\begin{equation}
0 \leq \sum_{i,j=1}^t M_{ij} p(F_i,G)p(F_j,G) = \sum_{i,j=1}^t M_{ij} p(F_i,F_j;G) + o(1).\label{eq:sigmaX}
\end{equation}
The inequality \eqref{eq:sigmaX} is usually applied to a large graph $G$ with a randomly chosen labeled vertices in a way that we now describe.
Fix $\sigma$-flags $F$ and $F'$ and a graph $G$.
We now define a random variable $p(F,F';G^\sigma)$ as follows:
label $|\sigma|$ vertices of $G$ with $1,\ldots,|\sigma|$ and if the resulting graph $G'$ is a $\sigma$-flag,
then $p(F',F';G^\sigma)=p(F,F';G')$;
if $G'$ is not a $\sigma$-flag, then $p(F_i,F_j;G^\sigma)=0$.
The expected value of $p(F,F';G^\sigma)$ can be expressed as a linear combination of
densities of $(|F|+|F'|-|\sigma|)$-vertex subgraphs of $G$~\cite{bib-razborov07},
i.e.,
there exist coefficients $\alpha_H$, $H\in\FF_{|F|+|F'|-|\sigma|}$, such that
\begin{equation}
\EE p(F,F';G^\sigma)=\sum_{H\in\FF_{|F|+|F'|-|\sigma|}}\alpha_H\cdot p(H,G)\label{eq:unlabX}
\end{equation}
for every graph $G$. It can be shown that $\alpha_H=\EE p(F,F';H^\sigma)$.

Let $\vec F = [F_1,\ldots,F_t]$ be a vector of $\ell$-vertex $\sigma$-flags and
let $M$ be a $t\times t$ positive semidefinite matrix.
The equality \eqref{eq:unlabX} yields that there exist coefficients $\alpha_H$ such that
\begin{equation}
\EE \sum_{i,j=1}^t M_{ij} p(F_i,F_j;G^\sigma)=\sum_{H\in\FF_{2\ell-|\sigma|}} \alpha_H\cdot p(H,G)\label{eq:almostX}
\end{equation}
for every graph $G$, which combines with \eqref{eq:sigmaX} to
\begin{equation}
0\le \sum_{H\in\FF_{2\ell-|\sigma|}} \alpha_H\cdot p(H,G)+o(1)\label{eq:finalX}
\end{equation}
for every graph $G$,
where
\[
\alpha_H = \sum_{i,j=1}^t M_{ij}\cdot\EE p(F_i,F_j;H^\sigma)
\]
In particular, the coefficients $\alpha_H$ depend only on the choice of $\vec F$ and $M$.

\section{Main result}

We start with proving the following lemma using the flag algebra method.

\begin{lemma}
\label{lm-fract}
Let $G$ be a weighted graph with all edges of weight one.
It holds that
$$\EEW\pi_{3,f}(G[W])\le 21+o(1)$$
where $W$ is a uniformly chosen random subset of seven vertices of $G$.
\end{lemma}

\begin{proof}
We use the flag algebra method to find coefficients $c_U$, $U\in\FF_7$, such that 
\begin{equation}
0 \leq \sum_{U \in \FF_7} c_U \cdot p(U,G) + o(1)\label{eq:c}
\end{equation}
and
\begin{equation}
\pi_{3,f}(U) + c_U \leq 21\label{eq:sharp}
\end{equation}
for every $U \in \FF_7$.
The statement of the lemma would then follow from \eqref{eq:c} and \eqref{eq:sharp} using $\sum_{U\in \FF_7} p(U,G) = 1$
as we next show.
\begin{align*}
\EEW\pi_{3,f}(G[W])
&=
\sum_{U \in \FF_7} \pi_{3,f}(U)\cdot p(U,G) \\
&\leq
\sum_{U \in \FF_7} (\pi_{3,f}(U)+c_U)\cdot p(U,G) + o(1)\\
&\leq \sum_{U \in \FF_7} 21\cdot p(U,G) + o(1) = 21 + o(1).
\end{align*}

We now focus on finding the coefficients $c_U$, $U\in\FF_7$, satisfying \eqref{eq:c} and \eqref{eq:sharp}.
Let $\sigma_1$ be a flag consisting of a single vertex labeled with $1$ and
consider the following vector $\vec F=(F_1,\ldots,F_7)$ of $\sigma_1$-flags from $\FF^{\sigma_1}_4$ (the single labeled
vertex is depicted by a white square and the remaining vertices by black circles).

\newcommand{\vc}[1]{\ensuremath{\vcenter{\hbox{#1}}}}
\tikzset{unlabeled_vertex/.style={inner sep=1.9pt, outer sep=0pt, circle, fill}} 
\tikzset{labeled_vertex/.style={inner sep=2.2pt, outer sep=0pt, rectangle, fill=white,draw}} 
\tikzset{edge_color1/.style={color=red,  line width=1.2pt,opacity=0}} 
\tikzset{edge_color2/.style={color=black, line width=0.8pt,opacity=1}} 
\def\outercycle#1#2{ \draw \foreach \x in {0,1,...,#1}{(270-45+\x*360/#2:0.7) coordinate(x\x)};}

\begin{align*}
\vec F = \left(
\vc{
{\begin{tikzpicture}\outercycle{5}{4}
\draw[edge_color1] (x0)--(x1);\draw[edge_color1] (x0)--(x2);\draw[edge_color1] (x0)--(x3);  \draw[edge_color1] (x1)--(x2);\draw[edge_color1] (x1)--(x3);  \draw[edge_color1] (x2)--(x3);    
\draw (x0) node[labeled_vertex]{};\draw (x1) node[unlabeled_vertex]{};\draw (x2) node[unlabeled_vertex]{};\draw (x3) node[unlabeled_vertex]{};
\end{tikzpicture}}
,
{\begin{tikzpicture}\outercycle{5}{4}
\draw[edge_color1] (x0)--(x1);\draw[edge_color1] (x0)--(x2);\draw[edge_color1] (x0)--(x3);  \draw[edge_color1] (x1)--(x2);\draw[edge_color1] (x1)--(x3);  \draw[edge_color2] (x2)--(x3);    
\draw (x0) node[labeled_vertex]{};\draw (x1) node[unlabeled_vertex]{};\draw (x2) node[unlabeled_vertex]{};\draw (x3) node[unlabeled_vertex]{};
\end{tikzpicture}}
,
{\begin{tikzpicture}\outercycle{5}{4}
\draw[edge_color1] (x0)--(x1);\draw[edge_color1] (x0)--(x2);\draw[edge_color2] (x0)--(x3);  \draw[edge_color1] (x1)--(x2);\draw[edge_color2] (x1)--(x3);  \draw[edge_color2] (x2)--(x3);    
\draw (x0) node[labeled_vertex]{};\draw (x1) node[unlabeled_vertex]{};\draw (x2) node[unlabeled_vertex]{};\draw (x3) node[unlabeled_vertex]{};
\end{tikzpicture}}
,
{\begin{tikzpicture}\outercycle{5}{4}
\draw[edge_color2] (x0)--(x1);\draw[edge_color2] (x0)--(x2);\draw[edge_color2] (x0)--(x3);  \draw[edge_color1] (x1)--(x2);\draw[edge_color1] (x1)--(x3);  \draw[edge_color1] (x2)--(x3);    
\draw (x0) node[labeled_vertex]{};\draw (x1) node[unlabeled_vertex]{};\draw (x2) node[unlabeled_vertex]{};\draw (x3) node[unlabeled_vertex]{};
\end{tikzpicture}}
,
{\begin{tikzpicture}\outercycle{5}{4}
\draw[edge_color1] (x0)--(x1);\draw[edge_color2] (x0)--(x2);\draw[edge_color2] (x0)--(x3);  \draw[edge_color1] (x1)--(x2);\draw[edge_color2] (x1)--(x3);  \draw[edge_color2] (x2)--(x3);    
\draw (x0) node[labeled_vertex]{};\draw (x1) node[unlabeled_vertex]{};\draw (x2) node[unlabeled_vertex]{};\draw (x3) node[unlabeled_vertex]{};
\end{tikzpicture}}
,
{\begin{tikzpicture}\outercycle{5}{4}
\draw[edge_color1] (x0)--(x1);\draw[edge_color2] (x0)--(x2);\draw[edge_color2] (x0)--(x3);  \draw[edge_color2] (x1)--(x2);\draw[edge_color2] (x1)--(x3);  \draw[edge_color1] (x2)--(x3);    
\draw (x0) node[labeled_vertex]{};\draw (x1) node[unlabeled_vertex]{};\draw (x2) node[unlabeled_vertex]{};\draw (x3) node[unlabeled_vertex]{};
\end{tikzpicture}}
,
{\begin{tikzpicture}\outercycle{5}{4}
\draw[edge_color1] (x0)--(x1);\draw[edge_color2] (x0)--(x2);\draw[edge_color2] (x0)--(x3);  \draw[edge_color2] (x1)--(x2);\draw[edge_color2] (x1)--(x3);  \draw[edge_color2] (x2)--(x3);    
\draw (x0) node[labeled_vertex]{};\draw (x1) node[unlabeled_vertex]{};\draw (x2) node[unlabeled_vertex]{};\draw (x3) node[unlabeled_vertex]{};
\end{tikzpicture}}
}
\right)
\end{align*}

Let $M$ be the following $7\times 7$-matrix.
\begin{center}
$M= \frac{1}{12\cdot 10^{9}}$
\scalebox{0.6}{
$\begin{pmatrix}
 1800000000 & 2444365956 &  640188285 & -1524146769 & 1386815580 & -732139362  & -129387078 \\
 2444365956 & 4759879134 & 1177441152 & -1783771230 & 2546923788 & -1397639394 & -143552208 \\
  640188285  & 1177441152 &  484273772 & -317303211&  1038156300  & -591902130   & -6783162 \\
-1524146769 & -1783771230 & -317303211 & 1558870290 & -651906630  & 305728704  & 154602378 \\
 1386815580 & 2546923788 & 1038156300 & -651906630 &  2285399634 & -1283125950  & -10755036 \\
 -732139362 & -1397639394 & -591902130 &  305728704 & -1283125950  & 734039016   & -1621938 \\
 -129387078  & -143552208  &  -6783162  & 154602378  & -10755036   & -1621938  &  23860164  
\end{pmatrix}$}\mbox{.}
\end{center}
The matrix $M$ is a positive semidefinite matrix with rank six;
the eigenvector corresponding to the zero eigenvalue is $(1, 0, 3, 1, 0, 3, 0)$.
Let
\[
c_U=\sum_{i,j=1}^7M_{ij}\EE p(F_i,F_j;U^{\sigma_1})\;\mbox{.}
\]
The inequality \eqref{eq:finalX} implies that
\[
0 \leq  \sum_{U \in \FF_7} c_U \cdot p(U,G) + o(1),
\]
which establishes~\eqref{eq:c}.
The inequality \eqref{eq:sharp} is verified with computer assistance
by evaluating the coefficient $c_U$ and the quantity $\pi_{3,f}(U)$ for each $U\in\FF_7$.
Since $|\FF_7|=1044$, we do not list $c_U$ and $\pi_{3,f}(U)$ here.
The computer programs that we used and their outputs have been made
available on arXiv as ancillary files and
are also available at \url{http://orion.math.iastate.edu/lidicky/pub/tile23}.
\end{proof}

The following lemma can be derived from the result of Haxell and R\"odl~\cite{bib-haxell01+} on fractional triangle decompositions or
from a more general result of Yuster~\cite{bib-yuster05}.

\begin{lemma}
\label{lm-decomp}
Let $G$ be a graph with $n$ vertices.
It holds that $\pi_3(G)\le\pi_{3,f}(G)+o(n^2)$.
\end{lemma}

We now use Lemmas~\ref{lm-fract} and~\ref{lm-decomp} to prove our main result.

\begin{theorem}
\label{thm-main}
Every $n$-vertex graph $G$ satisfies $\pi_3(G)\le (1/2+o(1))n^2$.
\end{theorem}

\begin{proof}
Fix an $n$-vertex graph $G$.
By Lemma~\ref{lm-decomp}, it is enough to show that $\pi_{3,f}(G)\le (1/2+o(1))n^2$.

Fix an optimal fractional $3$-de\-compo\-si\-tion of $G[W]$ for every $7$-vertex subset $W\subseteq V(G)$, and
set the weight $w(e)$ of an edge $e$ to the sum of its weights in the optimal fractional $3$-decomposition of $G[W]$
with $e\subseteq W$ multiplied by $\binom{n-2}{5}^{-1}$, and
the weight $w(t)$ of a triangle $t$ to the sum its weights in the optimal fractional $3$-decomposition of $G[W]$
with $t\subseteq W$ also multiplied by $\binom{n-2}{5}^{-1}$.
Since each edge $e$ of $G$ is contained in $\binom{n-2}{5}$ subsets $W$, we have obtained a fractional $3$-decomposition of $G$.
The weight of this decomposition is equal to
$$\frac{1}{\binom{n-2}{5}}\sum_{W\in\binom{V(G)}{7}}\pi_{3,f}(G[W])\le
  \frac{\binom{n}{7}}{\binom{n-2}{5}}(21 + o(1)) = n^2/2 + o(n^2)\;\mbox{,}$$
where the inequality follows from Lemma~\ref{lm-fract}.
We conclude that $\pi_{3,f}(G)\le n^2/2 + o(n^2)$, which completes the proof.
\end{proof}

The next corollary follows directly from Theorem~\ref{thm-main}.

\begin{corollary}
\label{cor-main}
Every $n$-vertex graph with $n^2/4+k$ edges contains $2k/3-o(n^2)$ edge-disjoint triangles.
\end{corollary}

\section{Concluding remarks}

Our first proof of this result, which can be found in~\cite{bib-orig-arxiv},
combined the flag algebra method and regularity method arguments.
In particular, we proved the fractional relaxation of Conjecture~\ref{conj}
in the setting of weighted graphs and with an additional restriction on its support;
this statement was then combined with a blow-up lemma for edge-decompositions recently proved
by Kim, K\"uhn, Osthus and Tyomkyn~\cite{bib-kim16+}.
It was then brought to our attention that the results from~\cite{bib-haxell01+,bib-yuster05}
allow obtaining our main result directly from the fractional relaxation,
which is the proof that we present here.
We believe that the argument using combinatorial designs that we applied in~\cite{bib-orig-arxiv}
to combine the flag algebra method and the blow-up lemma of Kim et al.~\cite{bib-kim16+}
can be of independent interest and so we wanted to mention the original proof of our result and
its idea here.

We also tried to prove Lemma~\ref{lm-fract} in the non-fractional setting,
i.e., to show that $\EEW\pi_{3}(G[W])\le 21+o(1)$.
Unfortunately, the computation with $7$-vertex flags yields only that $\EEW\pi_{3}(G[W])\le 21.588+o(1)$.
We would like to remark that if it were possible to prove Lemma~\ref{lm-fract} in the non-fractional setting,
we would be able to prove Theorem~\ref{thm-main} without using additional results as a blackbox:
we would consider a random $(n,7,2,1)$-design on the vertex set of an $n$-vertex graph $G$ as in~\cite{bib-orig-arxiv} and
apply the non-fractional version of Lemma~\ref{lm-fract} to this design.

Finally, we would also like to mention two open problems related to our main result.
Theorem~\ref{thm-main} asserts that $\pi_3(G)\le n^2/2+o(n^2)$ for every $n$-vertex graph $G$.
However, it could be true (cf.~the remark after Problem 41 in~\cite{bib-tuza91}) that
$\pi_3(G)\le n^2/2+2$ for every $n$-vertex graph $G$.
The second problem that we would like to mention is a possible generalization of Corollary~\ref{cor-main},
which is stated in~\cite{bib-tuza91} as Problem~42.
Fix $r\ge 4$. Does every $n$-vertex graph with $\frac{r-2}{2r-2}n^2+k$ edges
contain $\frac{2}{r}k-o(n^2)$ edge-disjoint complete graphs of order $r$?

\section*{Acknowledgements}

The authors would like to thank Ervin Gy\H{o}ri and Katherine Staden for their comments
on the problems considered in this paper.
The authors would also like to thank Allan Lo for drawing their attention to the paper~\cite{bib-haxell01+}.

\bibliographystyle{abbrv}
\bibliography{tile23}

\begin{thebibliography}{10}

\bibitem{bib-baber11}
R.~Baber and J.~Talbot.
\newblock Hypergraphs do jump.
\newblock {\em Combin. Probab. Comput.}, 20:161--171, 2011.

\bibitem{bib-baber14}
R.~Baber and J.~Talbot.
\newblock A solution to the 2/3 conjecture.
\newblock {\em SIAM J. Discrete Math.}, 28:756--766, 2014.

\bibitem{bib-balogh14}
J.~Balogh, P.~Hu, B.~Lidick\'y, and H.~Liu.
\newblock Upper bounds on the size of 4- and 6-cycle-free subgraphs of the
  hypercube.
\newblock {\em European J. Combin.}, 35:75--85, 2014.

\bibitem{bib-balogh16}
J.~Balogh, P.~Hu, B.~Lidick\'y, and F.~Pfender.
\newblock Maximum density of induced 5-cycle is achieved by an iterated blow-up
  of 5-cycle.
\newblock {\em European J. Combin.}, 52:47--58, 2016.

\bibitem{bib-balogh17}
J.~Balogh, P.~Hu, B.~Lidick\'y, F.~Pfender, J.~Volec, and M.~Young.
\newblock Rainbow triangles in three-colored graphs.
\newblock {\em J. Combin. Theory Ser. B}, 126:83--113, 2017.

\bibitem{bib-balogh15}
J.~Balogh, P.~Hu, B.~Lidick\'y, O.~Pikhurko, B.~Udvari, and J.~Volec.
\newblock Minimum number of monotone subsequences of length 4 in permutations.
\newblock {\em Combin. Probab. Comput.}, 24:658--679, 2015.

\bibitem{bib-chung81}
F.~R.~K. Chung.
\newblock On the decomposition of graphs.
\newblock {\em SIAM J. Algebraic Discrete Methods}, 2:1--12, 1981.

\bibitem{bib-coregliano17}
L.~N. Coregliano and A.~A. Razborov.
\newblock On the density of transitive tournaments.
\newblock {\em J. Graph Theory}, 85:12--21, 2017.

\bibitem{bib-cumm13}
J.~Cummings, D.~Kr\'al', F.~Pfender, K.~Sperfeld, A.~Treglown, and M.~Young.
\newblock Monochromatic triangles in three-coloured graphs.
\newblock {\em J. Combin. Theory Ser. B}, 103:489--503, 2013.

\bibitem{bib-das13}
S.~Das, H.~Huang, J.~Ma, H.~Naves, and B.~Sudakov.
\newblock A problem of {E}rd{\H o}s on the minimum number of {$k$}-cliques.
\newblock {\em J. Combin. Theory Ser. B}, 103:344--373, 2013.

\bibitem{bib-erdos66+}
P.~Erd\H{o}s, A.~W. Goodman, and L.~P\'osa.
\newblock The representation of a graph by set intersections.
\newblock {\em Canad. J. Math.}, 18:106--112, 1966.

\bibitem{bib-evan15}
C.~Even-Zohar and N.~Linial.
\newblock A note on the inducibility of 4-vertex graphs.
\newblock {\em Graphs Combin.}, 31:1367--1380, 2015.

\bibitem{bib-falgas15}
V.~Falgas-Ravry, E.~Marchant, O.~Pikhurko, and E.~R. Vaughan.
\newblock The codegree threshold for 3-graphs with independent neighborhoods.
\newblock {\em SIAM J. Discrete Math.}, 29:1504--1539, 2015.

\bibitem{bib-gethner17}
E.~Gethner, L.~Hogben, B.~Lidick\'y, F.~Pfender, A.~Ruiz, and M.~Young.
\newblock On crossing numbers of complete tripartite and balanced complete
  multipartite graphs.
\newblock {\em J. Graph Theory}, 84:552--565, 2017.

\bibitem{bib-glebov16}
R.~Glebov, D.~Kr\'al', and J.~Volec.
\newblock A problem of {E}rd{\H o}s and {S}\'os on 3-graphs.
\newblock {\em Israel J. Math.}, 211:349--366, 2016.

\bibitem{bib-goaoc15}
X.~Goaoc, A.~Hubard, R.~de~Joannis~de Verclos, J.-S. Sereni, and J.~Volec.
\newblock Limits of order types.
\newblock In {\em 31st {I}nternational {S}ymposium on {C}omputational
  {G}eometry}, volume~34 of {\em LIPIcs. Leibniz Int. Proc. Inform.}, pages
  300--314. Schloss Dagstuhl. Leibniz-Zent. Inform., Wadern, 2015.

\bibitem{bib-gyori15+}
E.~Gy\H{o}ri and B.~Keszegh.
\newblock On the number of edge-disjoint triangles in {$K_4$}-free graphs.
\newblock Preprint available as arXiv:1506.03306, 2015.

\bibitem{bib-gyori17+}
E.~Gy\H{o}ri and B.~Keszegh.
\newblock On the number of edge-disjoint triangles in {$K_4$}-free graphs.
\newblock {\em Electronic Notes in Discrete Mathematics}, 61:557--560, 2017.

\bibitem{bib-gyori79+}
E.~Gy\H{o}ri and A.~V. Kostochka.
\newblock On a problem of {G}. {O}. {H}. {K}atona and {T}. {T}arj\'an.
\newblock {\em Acta Math. Acad. Sci. Hungar.}, 34:321--327 (1980), 1979.

\bibitem{bib-gyori87+}
E.~Gy\H{o}ri and Z.~Tuza.
\newblock Decompositions of graphs into complete subgraphs of given order.
\newblock {\em Studia Sci. Math. Hungar.}, 22:315--320, 1987.

\bibitem{bib-hatami12}
H.~Hatami, J.~Hladk\'y, D.~Kr\'al', S.~Norine, and A.~Razborov.
\newblock Non-three-colourable common graphs exist.
\newblock {\em Combin. Probab. Comput.}, 21:734--742, 2012.

\bibitem{bib-haxell01+}
P.~E. Haxell and V.~R{\"{o}}dl.
\newblock Integer and fractional packings in dense graphs.
\newblock {\em Combinatorica}, 21:13--38, 2001.

\bibitem{bib-hladky17}
J.~Hladk\'y, D.~Kr\'al', and S.~Norin.
\newblock Counting flags in triangle-free digraphs.
\newblock {\em Combinatorica}, 37:49--76, 2017.

\bibitem{bib-kahn81}
J.~Kahn.
\newblock Proof of a conjecture of {K}atona and {T}arj\'an.
\newblock {\em Period. Math. Hungar.}, 12:81--82, 1981.

\bibitem{bib-kim16+}
J.~Kim, D.~K\"{u}hn, D.~Osthus, and M.~Tyomkyn.
\newblock A blow-up lemma for approximate decompositions.
\newblock Preprint available as arXiv:1604.07282, 2016.

\bibitem{bib-orig-arxiv}
D.~{Kr{\'a}l'}, B.~{Lidick{\'y}}, T.~L. {Martins}, and Y.~{Pehova}.
\newblock {Decomposing graphs into edges and triangles}.
\newblock {\em arXiv 1710.08486v2}, Oct. 2017.

\bibitem{bib-kral13}
D.~Kr\'al', C.-H. Liu, J.-S. Sereni, P.~Whalen, and Z.~B. Yilma.
\newblock A new bound for the {$2/3$} conjecture.
\newblock {\em Combin. Probab. Comput.}, 22:384--393, 2013.

\bibitem{bib-kral12}
D.~Kr\'al', L.~Mach, and J.-S. Sereni.
\newblock A new lower bound based on {G}romov's method of selecting heavily
  covered points.
\newblock {\em Discrete Comput. Geom.}, 48:487--498, 2012.

\bibitem{bib-kral14}
D.~Kr\'al' and O.~Pikhurko.
\newblock Quasirandom permutations are characterized by 4-point densities.
\newblock {\em Geom. Funct. Anal.}, 23:570--579, 2013.

\bibitem{bib-lidicky17+}
B.~Lidick\'y and F.~Pfender.
\newblock Semidefinite programming and {R}amsey numbers.
\newblock Preprint available as arXiv:1704.03592, 2017.

\bibitem{bib-mcguinness94b}
S.~McGuinness.
\newblock The greedy clique decomposition of a graph.
\newblock {\em J. Graph Theory}, 18:427--430, 1994.

\bibitem{bib-mcguinness94a}
S.~McGuinness.
\newblock Greedy maximum-clique decompositions.
\newblock {\em Combinatorica}, 14:335--343, 1994.

\bibitem{bib-razborov07}
A.~A. Razborov.
\newblock Flag algebras.
\newblock {\em J. Symbolic Logic}, 72:1239--1282, 2007.

\bibitem{bib-tuza91}
Z.~Tuza.
\newblock {Unsolved Combinatorial Problems, Part I}.
\newblock BRICS Lecture Series LS-01-1, 2001.

\bibitem{bib-yuster05}
R.~Yuster.
\newblock Integer and fractional packing of families of graphs.
\newblock {\em Random Structures \& Algorithms}, 26:110--118, 2005.

\end{thebibliography}


\end{document}